\documentclass[12pt,twoside]{article}
\usepackage{amssymb}
\font\teneufm=eufm10 scaled \magstep1
\font\seveneufm=eufm7 scaled \magstep1
\font\fiveeufm=eufm5  scaled \magstep1
\newfam\eufmfam
\textfont\eufmfam=\teneufm
\scriptfont\eufmfam=\seveneufm
\scriptscriptfont\eufmfam=\fiveeufm

\usepackage{mathrsfs}

\newfam\msbfam
\font\tenmsb=msbm10 scaled \magstep1  \textfont\msbfam=\tenmsb
\font\sevenmsb=msbm7 scaled \magstep1 \scriptfont\msbfam=\sevenmsb
\font\fivemsb=msbm5 scaled \magstep1  \scriptscriptfont\msbfam=\fivemsb

\def\RR{{\mathbb R}}
\def\CC{{\mathbb C}}

\def\PP{{\mathbb P}}

\def\ra{\rightarrow}

 \def\HollowBoxx #1#2#3{{\dimen0=#1 \advance\dimen0 by -#2
       \dimen1=#1 \advance\dimen1 by #3
        \vrule height 0pt depth #3 width #2
       \hskip -#3
       \vrule height #1 depth #3 width #3}}
 \def\LeftContraction{\mathord{\kern1.45pt \HollowBoxx{6pt}{3.5pt}{.4pt}}\,}

 \def\HollowBox #1#2#3{{\dimen0=#1 \advance\dimen0 by -#3
       \dimen1=#1 \advance\dimen1 by #3
        \vrule height #1 depth #3 width #3
        \vrule height 0pt depth #3 width #2
        \hskip -#3}}
 \def\RightContraction{\mathord{\, \HollowBox{6pt}{3.1pt}{.4pt}} \kern1.6pt}

\def\qed{{\hfill $\Box$}}
\newtheorem{theorem}{THEOREM}[section]

\newtheorem{lemma}[theorem]{Lemma}

\newtheorem{remark}[theorem]{Remark}

\begin{document}

\begin{center}
{\Large \bf On the Number of Affine Equivalence Classes
\vspace{0.3cm}\\
of Spherical Tube Hypersurfaces}
\medskip \\
\normalsize A. V. Isaev
\end{center}

\begin{quotation} \small \sl We consider Levi non-degenerate tube hypersurfaces in $\CC^{n+1}$ that are $(k,n-k)$-spherical, i.e. locally CR-equivalent to the hyperquadric with Levi form of signature $(k,n-k)$, with $n\le 2k$. We show that the number of affine equivalence classes of such hypersurfaces is infinite (in fact, uncountable) in the following cases: (i) $k=n-2$, $n\ge 7$;\linebreak (ii) $k=n-3$, $n\ge 7$; (iii) $k\le n-4$. For all other values of $k$ and $n$, except for $k=3$, $n=6$, the number of affine classes is known to be finite. The exceptional case $k=3$, $n=6$ has been recently resolved by Fels and Kaup who gave an example of a family of $(3,3)$-spherical tube hypersurfaces that contains uncountably many pairwise affinely non-equivalent elements. In this paper we deal with the Fels-Kaup example by different methods. We give a direct proof of the sphericity of the hypersurfaces in the Fels-Kaup family, and use the $j$-invariant to show that this family indeed contains an uncountable subfamily of pairwise affinely non-equivalent hypersurfaces.
\end{quotation}

\thispagestyle{empty}

\pagestyle{myheadings}
\markboth{A. V. Isaev}{Affine Equivalence Classes of Spherical Tube Hypersurfaces}

\setcounter{section}{0}

\section{Introduction}
\setcounter{equation}{0}

We consider real hypersurfaces in the complex space $\CC^{n+1}$ for $n\ge 1$. Specifically, we will be discussing {\it tube}\, hypersurfaces. Recall that a hypersurface $M$ is called tube if it has the form
$$
B+i\RR^{n+1},
$$
where $B$ is a hypersurface in $\RR^{n+1}\subset\CC^{n+1}$ called the {\it base}\, of $M$. The geometry of $M$ is completely determined by that of its base. Everywhere below we assume that $M$ is connected and $C^{\infty}$-smooth (in fact, for all our arguments smoothness of class $C^7$ is sufficient). 

Clearly, affine maps of the form
$$
Z\mapsto AZ+a,\quad A\in GL_{n+1}(\RR),\quad a\in\CC^{n+1},
$$
where $Z=(z_0,z_1,\dots,z_n)$ is a point in $\CC^{n+1}$, preserve the class of tube hypersurfaces. Two tube hypersurfaces that can be mapped onto each other by a map of this kind are called {\it affinely equivalent}. It is not hard to construct an example of two CR-equivalent tube hypersurfaces that are not affinely equivalent. Indeed, let $x_j:=\hbox{Re}\,z_j$, $j=0,\dots,n$, and define tube hypersurfaces $M_1$, $M_2$ for $n\ge 2$ as follows
\begin{equation}
\begin{array}{ll}
M_1: & \displaystyle x_0=\sum_{j=1}^{n-1}x_j^2-x_n^2,\\
\vspace{-0.3cm}\\
M_2: & \displaystyle x_0=\sum_{j=1}^{n-2}x_j^2+x_{n-1}x_n+x_n^3.
\end{array}\label{m1m2}
\end{equation}
One can show that these tube hypersurfaces are affinely non-equivalent. However, the polynomial transformation $\varphi:=\varphi_2^{-1}\circ\varphi_1$ maps $M_1$ onto $M_2$, where $\varphi_1$ and $\varphi_2$ are polynomial automorphisms of $\CC^{n+1}$ given, respectively, by the formulas
$$
\begin{array}{ll}
\varphi_1:&
\begin{array}{lll}
z_0^*&=&\displaystyle i\left(z_0-\frac{1}{2}\sum_{j=1}^{n-1}z_j^2+\frac{1}{2}z_n^2\right),\\
\vspace{-0.3cm}\\
z_j^*&=&\displaystyle\frac{z_j}{\sqrt{2}},\quad j=1,\dots,n,
\end{array}\\
\vspace{-0.3cm}\\
\hbox{and the formulas}&\\
\hspace{-0.3cm}\\
\varphi_2:&
\begin{array}{lll}
z_0^*&=&\displaystyle i\left(z_0-\frac{1}{2}\sum_{j=1}^{n-2}z_j^2-\frac{1}{2}z_{n-1}z_n-\frac{1}{4}z_n^3\right),\\
\vspace{-0.3cm}\\
z_j^*&=&\displaystyle\frac{z_j}{\sqrt{2}},\quad j=1,\dots,n-2,\\
\vspace{-0.3cm}\\
z_{n-1}^*&=&\displaystyle\frac{1}{\sqrt{2}}\left(z_n+\frac{1}{4}z_{n-1}+\frac{3}{8}z_n^2\right),\\
\vspace{-0.3cm}\\
z_n^*&=&\displaystyle\frac{1}{\sqrt{2}}\left(z_n-\frac{1}{4}z_{n-1}-\frac{3}{8}z_n^2\right)
\end{array}
\end{array}
$$
(here and everywhere below asterisks indicate the coordinates of the image of a point under a map). For a study of relationships between affine and CR-equivalence of tube hypersurfaces see \cite{Lo}. Related work on tube domains was done in \cite{KS1}, \cite{KS2}.

In fact, in the example above, the map $\varphi_j$ for each $j$ establishes CR-equivalence between $M_j$ and the quadric defined by the equation 
$$
\hbox{Im}\,z_0=\displaystyle\sum_{j=1}^{n-1}|z_j|^2-|z_n|^2.
$$
Tube hypersurfaces locally CR-equivalent to Levi non-degenerate quadrics are the subject of this paper. Let $M$ be a Levi non-degenerate hypersurface in $\CC^{n+1}$, and let its Levi form have $k$ positive and $n-k$ negative eigenvalues at every point, where we assume that $n\le 2k$. We then say that $M$ is {\it $(k,n-k)$-spherical}, if every point in $M$ has a neighborhood CR-equivalent to an open subset of the quadric
\begin{equation}
Q_{k,n-k}:=\left\{Z\in\CC^{n+1}:\displaystyle\hbox{Im}\,z_0=\langle z,z\rangle_{k,n-k}:=\sum_{j=1}^k|z_j|^2-\sum_{j=k+1}^n|z_j|^2\right\},\label{quadric}
\end{equation}
where $z:=(z_1,\dots,z_n)$. Recall that a $C^1$-smooth CR-diffeomorphism between $C^{\infty}$-smooth Levi non-degenerate hypersurfaces is in fact $C^{\infty}$-smooth (see e.g. \cite{Tu} and references therein); therefore, in the above definition the local CR-equivalence maps are assumed -- without loss of generality -- to be of class $C^{\infty}$.

In our terminology the hypersurfaces $M_1$, $M_2$ defined in (\ref{m1m2}) are $(n-1,1)$-spherical. Observe that $(n,0)$-spherical hypersurfaces are exactly hypersurfaces locally CR-equivalent to the sphere $S^{2n+1}\subset\CC^{n+1}$. We also remark that the quadric $Q_{k,n-k}$ for all $k$ and $n$ can be viewed as an algebraic quadratic tube hypersurface analogous to $M_1$ (see (\ref{m1m2})). Namely, the map
$$
\begin{array}{lll}
z_0^*&=&\displaystyle-iz_0+\sum_{j=1}^kz_j^2-\sum_{j=k+1}^nz_j^2,\\
\vspace{-0.3cm}\\
z^*&=&\sqrt{2}z,
\end{array}
$$
transforms $Q_{k,n-k}$ into the tube hypersurface given by the equation
$$
x_0=\displaystyle\sum_{j=1}^kx_j^2-\sum_{j=k+1}^nx_j^2.
$$

Spherical tube hypersurfaces possess remarkable properties. For example, as was shown in \cite{I4}, every such hypersurface is real-analytic and extends to a closed spherical tube hypersurface in $\CC^{n+1}$ (see also \cite{FK1}). Thus from now on we will only consider closed hypersurfaces. 

We are interested in the affine equivalence problem for spherical tube hypersurfaces. In \cite{Y} Yang proposed to approach this problem by means of utilizing the zero CR-curvature equations arising from the Cartan-Tanaka-Chern-Moser invariant theory (we note that an alternative approach has been recently proposed in \cite{FK2} -- see Remark \ref{admissalrrem}). In the tube case the zero CR-curvature equations significantly simplify and lead to a second-order system of partial differential equations for the defining function (see \cite{Y}, \cite{I3}, \cite{I5}). The properties of the solutions of this system can be thoroughly investigated in general (see \cite{I4}). For the strongly pseudoconvex case this approach had led to a complete explicit affine classification (see \cite{DY}). One consequence of this classification is the finiteness of the number of affine equivalence classes of $(n,0)$-spherical tube hypersurfaces for every $n\ge 1$.

Further, a complete affine classification of $(n-1,1)$-spherical tube hypersurfaces was obtained in \cite{IM} (hypersurfaces (\ref{m1m2}) are of course part of this classification). Again, for every $n\ge 2$ the number of affine equivalence classes turned out to be finite. We were also able to obtain a complete affine classification for the case $k=n-2$ (here $n\ge 4$). This classification was announced in paper \cite{I1}, where proofs were also briefly sketched. Full details were given in a very long preprint (see \cite{I2}). Because of the prohibitive length of the preprint complete proofs have never been published. 

Our classification of $(n-2,2)$-spherical tube hypersurfaces in \cite{I1}, \cite{I2} implies that the number of affine equivalence classes of such hypersurfaces is finite for $n\le 6$ and infinite for $n\ge 7$. One of the aims of the present paper is to provide a complete proof of the second part of this statement. We will see that algebraic non-quadratic hypersurfaces such as $M_2$ defined in (\ref{m1m2}) play an important role in the proof. In fact, our result is the following more general theorem (cf. Theorem 2 of \cite{I1}).

\begin{theorem}\label{theorem1}\sl The number of affine equivalence classes of closed connected $(k,n-k)$-spherical tube hypersurfaces in $\CC^{n+1}$ is infinite (in fact, uncountable) in the following cases: 
\vspace{0.1cm}

(i) $k=n-2$, $n\ge 7$; 

(ii) $k=n-3$, $n\ge 7$; 

(iii) $k\le n-4$.
\end{theorem}

\noindent The proof of Theorem \ref{theorem1} is given in Section \ref{prooftheorem1}. We explicitly present two real 1-parameter families of algebraic non-quadratic $(k,n-k)$-spherical tube hypersurfaces that cover cases (i)--(iii) of Theorem \ref{theorem1}, and show that all hypersurfaces in the families are pairwise affinely non-equivalent. Every hypersurface in each of the family is equivalent to the quadric $Q_{k,n-k}$ by means of a polynomial automorphism of $\CC^{n+1}$. We note that the techniques of \cite{FK2} yield an alternative proof of statements (ii) and (iii) of Theorem \ref{theorem1}, as well as an alternative way for obtaining the classifications of $(n,0)$- and $(n,1)$-spherical tube hypersurfaces found in \cite{DY}, \cite{IM} (see Remark \ref{admissalrrem}).

Theorem \ref{theorem1} does not cover the case $k=3$, $n=6$, and the question about the number of affine equivalence classes in this situation remained open until  recently. In \cite{FK2} Fels and Kaup constructed a real 1-parameter family of $(3,3)$-spherical hypersurfaces in $\CC^7$ and showed that it contains uncountably many pairwise affinely non-equivalent hypersurfaces. Let $S_t$ be the hypersurface in $\CC^7$ whose base is given by the equation
\begin{equation}
x_0=x_1x_6+x_2x_5+x_3x_4+x_4^3+x_5^3+x_6^3+tx_4x_5x_6,\quad t\in\RR.\label{hyperst}
\end{equation}
For every $t\in\RR$ consider the following cubic of three real variables
\begin{equation}
c_t(u_1,u_2,u_3):=u_1^3+u_2^3+u_3^3+tu_1u_2u_3.\label{cubics}
\end{equation}
To every cubic $c_t$ Fels and Kaup assign a certain associative commutative algebra ${\mathcal V}_{c_t}$ with the property that two such algebras ${\mathcal V}_{c_{t_1}}$, ${\mathcal V}_{c_{t_2}}$ are isomorphic if and only if the cubics $c_{t_1}$, $c_{t_2}$ are linearly equivalent (i.e. can be mapped into each other by a transformation from $GL_3(\RR)$). On the other hand, ${\mathcal V}_{c_{t_1}}$, ${\mathcal V}_{c_{t_2}}$ are isomorphic if and only if the hypersurfaces $S_{t_1}$, $S_{t_2}$ are affinely equivalent. It is shown in \cite{FK2} that the cubics $c_t$ are pairwise linearly non-equivalent, for example, for $t$ lying in a certain interval $I\subset\RR$. Furthermore, the $(3,3)$-sphericity of every hypersurface $S_t$ is a consequence of the general algebraic approach of \cite{FK1}, \cite{FK2} to classifying all local tube realizations of a given CR-manifold (quadric (\ref{quadric}) in his case). Thus one obtains the following result.

\begin{theorem}\label{theorem2}{\rm\cite{FK2}} \sl The hypersurfaces $S_t$ are $(3,3)$-spherical for all $t\in\RR$ and there exists an interval $I\subset\RR$ such that $S_t$ are pairwise affinely non-equivalent for $t\in I$. In particular, the number of affine equivalence classes of closed connected $(3,3)$-spherical tube hypersurfaces in $\CC^7$ is infinite (in fact, uncountable).
\end{theorem}

In this paper (see Section \ref{prooftheorem2}) we give a direct (and mostly analytic) proof of Theorem \ref{theorem2} different from that in \cite{FK2}. First of all, we show that every $S_t$ is $(3,3)$-spherical by explicitly presenting a polynomial automorphism of $\CC^7$ that transforms $S_t$ into $Q_{3,3}$. Next, we use the same general method as in the proof of Theorem \ref{theorem1} to show that linear non-equivalence over $\RR$ of two cubics $c_{t_1}$, $c_{t_2}$ implies affine non-equivalence of the corresponding hypersurfaces $S_{t_1}$, $S_{t_2}$. Finally, we establish pairwise linear non-equivalence of the cubics $c_t$ for small $|t|$ by comparing the values of the $j$-invariant for the elliptic curves defined as the zero loci of $c_t$ in $\CC\PP^2$.

{\bf Acknowledgements.} This work is supported by the Australian Research Council and part of it was done during the author's visit to the Indian Institute of Science in Bangalore. The author is grateful to M. G. Eastwood for many useful discussions and for suggesting, on a number of occasions, shorter and more elegant proofs.
 
\section{Proof of Theorem \ref{theorem1}}\label{prooftheorem1}
\setcounter{equation}{0}

Assume first that $5\le k\le n-2$, $n\ge 7$, and consider the family of algebraic hypersurfaces $P_t$ given by the following equation

\newpage

\begin{equation}
\begin{array}{l}
\displaystyle x_0=\sum_{j=1}^{k-2}x_j^2+x_{k-1}x_k+x_{k+1}x_{k+2}-\sum_{j=k+3}^nx_j^2+\\
\vspace{-0.1cm}\\
\displaystyle\hspace{1cm}2\sqrt{2(1+t)}x_{k-4}x_{k-1}x_{k+1}+2\sqrt{3t}\,x_{k-3}x_{k+1}^2+\frac{1+t}{\sqrt{3t}}x_{k-3}x_{k-1}^2+\\
\vspace{-0.1cm}\\
\displaystyle\hspace{1cm}\sqrt{\frac{-t^2+34t-1}{3t}}x_{k-2}x_{k-1}^2+(x_{k-1}^2+x_{k+1}^2)(x_{k-1}^2+tx_{k+1}^2),
\end{array}\label{surfht}
\end{equation}
where $1\le t\le 17+12\sqrt{2}$. Every hypersurface $P_t$ is $(k,n-k)$-spherical. Indeed, $P_t$ is mapped onto $Q_{k,n-k}$ by the following polynomial automorphism $p_t$ of $\CC^{n+1}$
$$
\begin{array}{l}
\displaystyle z_0^*=i\Biggl(z_0-\frac{1}{2}\sum_{j=1}^{k-2}z_j^2-\frac{1}{2}z_{k-1}z_k-\frac{1}{2}z_{k+1}z_{k+2}+\frac{1}{2}\sum_{j=k+3}^nz_j^2-\\
\vspace{-0.3cm}\\
\displaystyle\hspace{1cm}\frac{\sqrt{2(1+t)}}{2}z_{k-4}z_{k-1}z_{k+1}-\frac{\sqrt{3t}}{2}z_{k-3}z_{k+1}^2-\frac{1+t}{4\sqrt{3t}}z_{k-3}z_{k-1}^2-\\
\vspace{-0.3cm}\\
\displaystyle\hspace{1cm}\frac{1}{4}\sqrt{\frac{-t^2+34t-1}{3t}}z_{k-2}z_{k-1}^2-\frac{1}{8}(z_{k-1}^2+z_{k+1}^2)(z_{k-1}^2+tz_{k+1}^2)\Biggr),\\
\vspace{-0.1cm}\\
\displaystyle z_j^*=\frac{1}{\sqrt{2}}z_j,\quad j=1,\dots,k-5,k+3,\dots,n,\\
\vspace{-0.1cm}\\
\displaystyle z_{k-4}^*=\frac{1}{\sqrt{2}}\left(z_{k-4}+\frac{1}{2}\sqrt{2(1+t)}z_{k-1}z_{k+1}\right),\\
\vspace{-0.1cm}\\
\displaystyle z_{k-3}^*=\frac{1}{\sqrt{2}}\left(z_{k-3}+\frac{\sqrt{3t}}{2}z_{k+1}^2+\frac{1+t}{4\sqrt{3t}}z_{k-1}^2\right),\\
\vspace{-0.1cm}\\
\displaystyle z_{k-2}^*=\frac{1}{\sqrt{2}}\left(z_{k-2}+\frac{1}{4}\sqrt{\frac{-t^2+34t-1}{3t}}z_{k-1}^2\right),\\
\vspace{-0.1cm}\\\displaystyle z_{k-1}^*=-\frac{i}{4}\left(2z_{k-1}+z_k+\sqrt{\frac{-t^2+34t-1}{3t}}z_{k-2}z_{k-1}+\right.\\
\vspace{-0.3cm}\\
\displaystyle\hspace{1cm}\left.\sqrt{2(1+t)}z_{k-4}z_{k+1}+\frac{1+t}{\sqrt{3t}}z_{k-3}z_{k-1}+z_{k-1}^3+\frac{1+t}{2}z_{k-1}z_{k+1}^2\right),\\
\vspace{-0.1cm}\\\displaystyle z_k^*=-\frac{i}{4}\left(2z_{k+1}+z_{k+2}+\sqrt{2(1+t)}z_{k-4}z_{k-1}+2\sqrt{3t}\,z_{k-3}z_{k+1}+\right.\\
\vspace{-0.3cm}\\
\displaystyle\hspace{1cm}\left.\frac{1+t}{2}z_{k-1}^2z_{k+1}+tz_{k+1}^3\right),\\
\end{array}
$$
$$
\begin{array}{l}

\displaystyle z_{k+1}^*=-\frac{i}{4}\left(-2z_{k-1}+z_k+\sqrt{\frac{-t^2+34t-1}{3t}}z_{k-2}z_{k-1}+\right.\\
\vspace{-0.3cm}\\
\displaystyle\hspace{1cm}\left.\sqrt{2(1+t)}z_{k-4}z_{k+1}+\frac{1+t}{\sqrt{3t}}z_{k-3}z_{k-1}+z_{k-1}^3+\frac{1+t}{2}z_{k-1}z_{k+1}^2\right),\\
\vspace{-0.1cm}\\
\displaystyle z_{k+2}^*=-\frac{i}{4}\left(-2z_{k+1}+z_{k+2}+\sqrt{2(1+t)}z_{k-4}z_{k-1}+2\sqrt{3t}\,z_{k-3}z_{k+1}+\right.\\
\vspace{-0.3cm}\\
\displaystyle\hspace{1cm}\left.\frac{1+t}{2}z_{k-1}^2z_{k+1}+tz_{k+1}^3\right).
\end{array}
$$

The above formulas for the map $p_t$ are not a result of guessing. They arise when one attempts to reduce the equation of the hypersurface $P_t$ to the Chern-Moser normal form by performing transformations specified in Lemmas 3.2 and 3.3 of \cite{CM}. It turns out that these transformations are sufficient for normalizing the equation of $P_t$, and no further (harder) steps of the Chern-Moser normalization process are necessary. This is the case for all known examples of algebraic spherical tube hypersurfaces and may be a manifestation of a general fact.

Next, we will show that the hypersurfaces $P_{t_1}$ and $P_{t_2}$ are affinely non-equivalent for $t_1\ne t_2$. In order to see this, we will need the following general lemma which is also of independent interest. 

\begin{lemma}\label{affhom}\sl Any hypersurface in $\RR^{n+1}$ given by an equation of the form
\begin{equation}
\begin{array}{l}
\displaystyle x_0=\sum_{j=1}^{k-2}x_j^2+x_{k-1}x_k+x_{k+1}x_{k+2}-\sum_{j=k+3}^nx_j^2+\\
\vspace{-0.1cm}\\
\displaystyle\hspace{1cm}ax_{k-4}x_{k-1}x_{k+1}+b\,x_{k-3}x_{k+1}^2+cx_{k-3}x_{k-1}^2+\\
\vspace{0.1cm}\\
\displaystyle\hspace{1cm}dx_{k-2}x_{k-1}^2+{\mathcal Q}_4(x_{k-1},x_{k+1}),
\end{array}\label{genhyper}
\end{equation}
where $5\le k\le n-2$, ${\mathcal Q}_4$ is a homogeneous polynomial of degree 4, and $a,b,c\in\RR^*$, $d\in\RR$, is affinely homogeneous.
\end{lemma}

\noindent {\bf Proof:} It is sufficient to prove the lemma for $k=5$, $n=7$. Let $S$ be a hypersurface of the form (\ref{genhyper}) in $\RR^8$. We apply to $S$ the translation $\tau_q: X\mapsto X-q$, where $X:=(x_0,\dots,x_7)$ and $q\in S$. Absorbing linear terms into $x_0$, we turn the equation of $\tau_q(S)$ into an equation of the form 
$$
\begin{array}{l}
\displaystyle x_0=x_1^2+x_2^2+x_3^2+x_4x_5+x_6x_7+L_1(x_1,x_2,x_3,x_4)x_4+\\
\vspace{-0.3cm}\\
\hspace{1cm}L_2(x_1,x_2,x_4,x_6)x_6+ax_1x_4x_6+bx_2x_6^2+cx_2x_4^2+dx_3x_4^2+\\
\vspace{-0.3cm}\\
\hspace{1cm}Ax_4^3+Bx_6^3+Cx_4x_6^2+Dx_4^2x_6+{\mathcal Q}_4(x_4,x_6),
\end{array}
$$
where $L_1$ and $L_2$ are linear functions and $A,B,C,D\in\RR$. Absorbing $L_1$ into $x_5$ and $L_2$ into $x_7$, and replacing $x_2$ by $x_2-(A/c)\,x_4-(B/b)\,x_6$, we obtain the equation 
$$
\begin{array}{l}
\displaystyle x_0=x_1^2+x_2^2+x_3^2+x_4x_5+x_6x_7+L_3(x_2,x_4)x_4+L_4(x_2,x_4,x_6)x_6+\\
\vspace{-0.3cm}\\
\hspace{1cm} ax_1x_4x_6+bx_2x_6^2+cx_2x_4^2+dx_3x_4^2+C'x_4x_6^2+D'x_4^2x_6+{\mathcal Q}_4(x_4,x_6),
\end{array}
$$    
where $L_3$ and $L_4$ are linear functions and $C',D'\in\RR$. Further, replacing $x_1$ by $x_1-(D'/a)\,x_4-(C'/a)x_6$, we get
$$
\begin{array}{l}
\displaystyle x_0=x_1^2+x_2^2+x_3^2+x_4x_5+x_6x_7+L_5(x_1,x_2,x_4)x_4+L_6(x_1,x_2,x_4,x_6)x_6+\\
\vspace{-0.3cm}\\
\hspace{1cm} ax_1x_4x_6+bx_2x_6^2+cx_2x_4^2+dx_3x_4^2+{\mathcal Q}_4(x_4,x_6),
\end{array}
$$
where $L_5$ and $L_6$ are linear functions. Finally, absorbing $L_5$ into $x_5$ and $L_6$ into $x_7$, we obtain the original equation of $S$. This argument shows that $S$ is affinely homogeneous, as required.\qed
\vspace{0.5cm} 

Let $B_{t_1}$ and $B_{t_2}$ be the bases of $P_{t_1}$ and $P_{t_2}$, respectively. Lemma \ref{affhom} implies that the hypersurfaces $B_{t_1}$ and $B_{t_2}$ are affinely homogeneous, hence $B_{t_1}$ and $B_{t_2}$ are affinely equivalent if and only if they are linearly equivalent. It is immediate from (\ref{surfht}) that the cubic terms in the equations of $B_{t_1}$ and $B_{t_2}$ are trace-free, where the trace is calculated with respect to the non-degenerate quadratic form $\sum_{j=1}^{k-2}x_j^2+x_{k-1}x_k+x_{k+1}x_{k+2}-\sum_{j=k+3}^nx_j^2$. Therefore Proposition 1 of \cite{EE} (which is based on an argument contained in \cite{Le}) yields that any linear equivalence between $B_{t_1}$ and $B_{t_2}$ has the form
$$
x_0\mapsto\lambda x_0,\quad x\mapsto Cx,
$$
where $x:=(x_1,\dots,x_n)$, $\lambda\in\RR^*$, $C\in GL_n(\RR)$. It then follows that a sufficient condition for affine non-equivalence of $B_{t_1}$ and $B_{t_2}$ is linear non-equivalence over $\RR$ of the fourth-order terms in their equations. Thus to establish pairwise affine non-equivalence of the hypersurfaces $P_t$ it is sufficient to establish pairwise linear non-equivalence of the quartics $q_t(\xi,\eta):=(\xi^2+\eta^2)(\xi^2+t\eta^2)$ over $\RR$.

Suppose that there exists a non-degenerate linear map $\varphi$
$$
\xi^*=\alpha \xi+\beta \eta,\quad \eta^*=\gamma \xi+\delta \eta,\quad\alpha,\beta,\gamma,\delta\in\RR
$$
that transforms $q_{t_1}$ into $q_{t_2}$, with $t_1<t_2$. We allow $\xi$ and $\eta$ to be complex and consider $\varphi$ as a transformation of $\CC^2$. Then $\varphi$ maps the zero locus of $q_{t_1}$ into the zero locus of $q_{t_2}$ (for any $t$ the zero locus of $q_t$ consists of the complex lines $\{\xi=\pm i\eta\}$, $\{\xi=\pm i\sqrt{t}\eta\}$). Consider the M\"obius transformation $m_{\varphi}$ of $\CC\PP^1$ arising from $\varphi$. Clearly, on the subset $\eta\ne 0$ of $\CC\PP^1$ for $\zeta=\xi/\eta$ we have
$$
\displaystyle m_{\varphi}(\zeta)=\frac{\alpha\zeta+\beta}{\gamma\zeta+\delta},
$$     
and $m_{\varphi}$ maps the set $\sigma_{t_1}:=\{\pm i,\pm i\sqrt{t_1}\}$ onto the set $\sigma_{t_2}:=\{\pm i,\pm i\sqrt{t_2}\}$. 

If $t_1=1$, then $\sigma_{t_1}$ is a two-point set and cannot be mapped onto the four-point set $\sigma_{t_2}$, hence such a map $\varphi$ does not exist for $t_1=1$. Assume now that $t_1>1$. Then $m_{\varphi}$ preserves the imaginary axis in the $\zeta$-plane, which immediately implies that $\alpha\gamma=0$ and $\beta\delta=0$. Then we have either $\beta=\gamma=0$ or $\alpha=\delta=0$, thus $m_{\varphi}$ is either a real dilation or the composition of a real dilation and $1/\zeta$. In either case $m_{\varphi}$ cannot map $\sigma_{t_1}$ onto $\sigma_{t_2}$ since $t_1<t_2$. 

We have thus shown that $q_{t_1}$ and $q_{t_2}$ are not equivalent by means of a transformation from $GL_2(\RR)$. Note that this statement no longer holds if maps from $GL_2(\CC)$ are allowed.    

Thus the algebraic spherical tube hypersurfaces $P_t$ defined in (\ref{surfht}) are pairwise affinely non-equivalent. It then follows that the number of affine equivalence classes of $(k,n-k)$-spherical tube hypersurfaces is uncountable in the following situations: (a) $k=n-2$, $n\ge 7$; (b) $k=n-3$, $n\ge 8$; (c) $k=n-4$, $n\ge 9$; (d) $k\le n-5$.

It remains to prove the theorem for $k=4$, $n=7$ and for $k=4$, $n=8$. We will present another real 1-parameter family of hypersurfaces that is defined in a more general setting than these two remaining cases. Assume that $4\le k\le n-3$, $n\ge 7$ and consider the family of hypersurfaces ${\mathcal P}_t$ given by the equation
$$
\begin{array}{l}
\displaystyle x_0=\sum_{j=1}^{k-2}x_j^2-x_{k-1}^2+x_kx_{k+1}+x_{k+2}x_{k+3}-\sum_{j=k+4}^nx_j^2+\\
\vspace{-0.1cm}\\
\displaystyle\hspace{1cm}2\sqrt{2(1+t)}x_{k-3}x_kx_{k+2}+2\sqrt{3t}\,x_{k-2}x_{k+2}^2+\frac{1+t}{\sqrt{3t}}x_{k-2}x_k^2+\\
\vspace{-0.1cm}\\
\displaystyle\hspace{1cm}\sqrt{\frac{t^2-34t+1}{3t}}x_{k-1}x_k^2+(x_k^2+x_{k+2}^2)(x_k^2+tx_{k+2}^2),
\end{array}
$$
where $t\ge 17+12\sqrt{2}$. Every hypersurface ${\mathcal P}_t$ is $(k,n-k)$-spherical. Indeed, ${\mathcal P}_t$ is mapped onto $Q_{k,n-k}$ by the following polynomial automorphism of $\CC^{n+1}$
$$
\begin{array}{l}
\displaystyle z_0^*=i\Biggl(z_0-\frac{1}{2}\sum_{j=1}^{k-2}z_j^2+\frac{1}{2}z_{k-1}^2-\frac{1}{2}z_kz_{k+1}-\frac{1}{2}z_{k+2}z_{k+3}+\frac{1}{2}\sum_{j=k+4}^nz_j^2-\\
\vspace{-0.3cm}\\
\displaystyle\hspace{1cm}\frac{\sqrt{2(1+t)}}{2}z_{k-3}z_kz_{k+2}-\frac{\sqrt{3t}}{2}z_{k-2}z_{k+2}^2-\frac{1+t}{4\sqrt{3t}}z_{k-2}z_k^2-\\
\vspace{-0.3cm}\\
\displaystyle\hspace{1cm}\frac{1}{4}\sqrt{\frac{t^2-34t+1}{3t}}z_{k-1}z_k^2-\frac{1}{8}(z_k^2+z_{k+2}^2)(z_k^2+tz_{k+2}^2)\Biggr),\\
\end{array}
$$
$$
\begin{array}{l}
\displaystyle z_j^*=\frac{1}{\sqrt{2}}z_j,\quad j=1,\dots,k-4,\,k+4,\dots,n,\\
\vspace{-0.1cm}\\
\displaystyle z_{k-3}^*=\frac{1}{\sqrt{2}}\left(z_{k-3}+\frac{1}{2}\sqrt{2(1+t)}z_kz_{k+2}\right),\\
\vspace{-0.1cm}\\
\displaystyle z_{k-2}^*=\frac{1}{\sqrt{2}}\left(z_{k-2}+\frac{\sqrt{3t}}{2}z_{k+2}^2+\frac{1+t}{4\sqrt{3t}}z_k^2\right),\\
\vspace{-0.1cm}\\
\displaystyle z_{k-1}^*=-\frac{i}{4}\left(2z_{k+2}+z_{k+3}+\sqrt{2(1+t)}z_{k-3}z_k+2\sqrt{3t}\,z_{k-2}z_{k+2}+\right.\\
\vspace{-0.3cm}\\
\displaystyle\hspace{1cm}\left.\frac{1+t}{2}z_k^2z_{k+2}+tz_{k+2}^3\right),\\
\vspace{-0.3cm}\\
\displaystyle z_k^*=-\frac{i}{4}\left(2z_k+z_{k+1}+\sqrt{\frac{t^2-34t+1}{3t}}z_{k-1}z_k+\right.\\
\vspace{-0.3cm}\\
\displaystyle\hspace{1cm}\left.\sqrt{2(1+t)}z_{k-3}z_{k+2}+\frac{1+t}{\sqrt{3t}}z_{k-2}z_k+z_k^3+\frac{1+t}{2}z_kz_{k+2}^2\right),\\
\vspace{-0.1cm}\\
\displaystyle z_{k+1}^*=\frac{1}{\sqrt{2}}\left(z_{k-1}-\frac{1}{4}\sqrt{\frac{t^2-34t+1}{3t}}z_k^2\right),\\
\vspace{-0.1cm}\\
\displaystyle z_{k+2}^*=-\frac{i}{4}\left(-2z_k+z_{k+1}+\sqrt{\frac{t^2-34t+1}{3t}}z_{k-1}z_k+\right.\\
\vspace{-0.3cm}\\
\displaystyle\hspace{1cm}\left.\sqrt{2(1+t)}z_{k-3}z_{k+2}+\frac{1+t}{\sqrt{3t}}z_{k-2}z_k+z_k^3+\frac{1+t}{2}z_kz_{k+2}^2\right),\\
\vspace{-0.3cm}\\
\displaystyle z_{k+3}^*=-\frac{i}{4}\left(-2z_{k+2}+z_{k+3}+\sqrt{2(1+t)}z_{k-3}z_k+2\sqrt{3t}\,z_{k-2}z_{k+2}+\right.\\
\vspace{-0.3cm}\\
\displaystyle\hspace{1cm}\left.\frac{1+t}{2}z_k^2z_{k+2}+tz_{k+2}^3\right).
\end{array}
$$ 

A proof completely analogous to that for the hypersurfaces $P_t$ above shows that the hypersurfaces ${\mathcal P}_t$ are affinely homogeneous and pairwise affinely non-equivalent. It then follows that the number of affine equivalence classes of $(k,n-k)$-spherical tube hypersurfaces is uncountable in the cases $k=4$, $n=7$ and $k=4$, $n=8$ as well.

The proof of the theorem is complete.\qed

\begin{remark}\label{unification}\rm For $t\ne 17+12\sqrt{2}$ one can write the families of hypersurfaces $P_t$ and ${\mathcal P}_t$ in a unified form. Namely, consider the family ${\mathfrak P}_{\tau}^p$ of tube hypersurfaces in $\CC^n$, with $n\ge 7$, given by the equations
$$
\begin{array}{l}
\displaystyle x_0=4x_1x_7+4x_2x_6-\tau x_3^2+2x_4^2-\tau x_5^2+4x_3x_5+\\
\vspace{-0.3cm}\\
\displaystyle\hspace{1cm}\sum_{j=8}^{p+7}x_j^2-\sum_{j=p+8}^nx_j^2-
2\tau x_1^2x_3-2\tau x_2^2x_5+4x_1^2x_5+4x_2^2x_3+\\
\vspace{-0.1cm}\\
\displaystyle\hspace{1cm}8x_1x_2x_4-\frac{\tau}{3}x_1^4+4x_1^2x_2^2-\frac{\tau}{3}x_2^4,
\end{array}
$$
where $\tau\in\RR$, $\tau\ne\pm 2$, and $0\le p\le n-7$. Note that every hypersurface ${\mathfrak P}_{\tau}^p$ is affinely equivalent to a hypersurface ${\mathfrak P}_{\tau'}^p$ for which $\tau'\in[-6,-2)\cup(-2,2)$. The Levi form of ${\mathfrak P}_{\tau}^p$ has signature $(p+5,n-(5+p))$ for $\tau\in[-6,-2)$ and signature $(p+4,n-(p+4))$ for $\tau\in(-2,2)$. By using Lemmas 3.2 and 3.3 of \cite{CM} as in the proof of Theorem \ref{theorem1}, one can show that ${\mathfrak P}_{\tau}^p$ is $(p+5,n-(5+p))$-spherical for $\tau\in[-6,-2)$ and is $(p+4,n-(p+4))$-spherical for $\tau\in(-2,2)$. Next, observe that the map
$$
\chi: t\mapsto -\frac{12\sqrt{t}}{t+1}
$$
is a bijection from $[1,17+12\sqrt{2})$ onto $[-6,-2)$ and from $(17+12\sqrt{2},\infty)$ onto $(-2,0)$. It is now not hard to see that ${\mathfrak P}_{\tau}^{k-5}$ is affinely equivalent to $P_{\chi^{-1}(\tau)}$ for $\tau\in[-6,-2)$, $5\le k\le n-2$, and that ${\mathfrak P}_{\tau}^{k-4}$ is affinely equivalent to ${\mathcal P}_{\chi^{-1}(\tau)}$ for $\tau\in(-2,0)$, $4\le k\le n-3$. In particular, the family ${\mathfrak P}_{\tau}^{k-4}$ for $\tau\in(-2,2)$, $4\le k\le n-3$ extends the family ${\mathcal P}_t$.         
\end{remark}

\begin{remark}\label{admissalrrem}\rm It is shown in \cite{FK2} that the affine classification problem for spherical tube hypersurfaces reduces to the problem of classifying finite-dimensional real and complex nilpotent associative algebras with\linebreak 1-dimensional annihilator up to abstract algebra isomorphisms (for brevity let us call such algebras {\it admissible}). In particular, there is an algorithm that associates a spherical tube hypersurface to every admissible algebra, and hypersurfaces associated to two such algebras are affinely equivalent if and only if the algebras are isomorphic. Although a general classification of admissible algebras does not exist, partial classifications corresponding to the cases $k=n$ and $k=n-1$ are not hard to obtain. They allow to recover the affine classifications of $(n,0)$- and $(n-1,1)$-spherical hypersurfaces found \cite{DY}, \cite{IM} without resorting to the zero CR-curvature equations (see \cite{FK2}). 

Interestingly, admissible algebras arise from isolated hypersurface singularities. Suppose that $W$ is a complex hypersurface in $\CC^N$ that has an isolated singularity at the origin given as the zero locus of a holomorphic function $f(z_1,\dots,z_N)$. Consider the corresponding {\it moduli algebra}
$$
{\mathcal A}:={\mathcal O}_{\CC^N,0}/I(f),
$$
where ${\mathcal O}_{\CC^N,0}$ is the algebra of functions in $z_1,\dots,z_N$ holomorphic in a neighborhood of the origin, and $I(f)$ is the ideal in ${\mathcal O}_{\CC^N,0}$ generated by $f$ and all its partial derivatives. It is well-known that ${\mathcal A}$ is a finite-dimensional commutative associative algebra depending only on the germ of $W$ at the origin (see e.g. \cite{GLS}, Section 2.1). Nakayama's lemma then implies that the (unique) maximal ideal ${\mathcal I}$ of ${\mathcal A}$ is a complex nilpotent associative algebra. Assume now that the origin is a quasi-homogeneous singularity (see \cite{S1}). In this case $f$ lies in its Jacobian ideal $J(f)$ generated by the partial derivatives of $f$. Hence $I(f)=J(f)$, and Lemma 3.3 of \cite{S2} implies that ${\mathcal I}$ is a complex admissible algebra. Thus, one can construct spherical tube hypersurfaces from ${\mathcal I}$ and its real forms by applying the algorithm of \cite{FK2}. 

For example, for $n=7$ and every $\tau\in[-6,-2)\cup(-2,2)$ the hypersurface ${\mathfrak P}_{\tau}^0$ introduced in Remark \ref{unification} arises from a real form of the maximal ideal of the moduli algebra of a simple elliptic hypersurface singularity of type $\tilde E_7$. Similarly, for $t\ne 0,6$ the hypersurface $S_t$ defined in (\ref{hyperst}) arises from a real form of the maximal ideal of the moduli algebra of a simple elliptic hypersurface singularity of type $\tilde E_6$ (for details concerning simple elliptic singularities see \cite{S2}). This real form is isomorphic to the algebra ${\mathcal V}_{c_{-18/t}}$ mentioned in the introduction. It is interesting to observe that if, in addition, $t\ne -3$, then the values of the $j$-invariant of the elliptic curves defined by the zero loci of $c_t$ and $c_{-18/t}$ in $\CC\PP^2$ are reciprocal (see \cite{E}). The zero loci of the cubics $c_t$ for $t\in\RR$, $t\ne -3$, as well as the $j$-invariant, play a role in our proof of Theorem \ref{theorem2} in the next section.
\end{remark}

\section{Proof of Theorem \ref{theorem2}}\label{prooftheorem2}
\setcounter{equation}{0}

In this section it will be convenient for us to replace the Hermitian form $\langle z,z\rangle_{3,3}$ by the equivalent Hermitian form
$$
\langle z,z\rangle_{3,3}':=\frac{1}{2}\,\hbox{Re}\,\left(z_1\overline{z}_6+z_2\overline{z}_5+z_3\overline{z}_4\right),
$$
and the quadric $Q_{3,3}$ by the equivalent quadric
$$
Q_{3,3}':=\left\{Z\in\CC^7:\hbox{Im}\,z_0=\langle z,z\rangle_{3,3}'\right\}.
$$

First of all, we observe that every hypersurface $S_t$ defined in (\ref{hyperst}) is $(3,3)$-spherical. Indeed, $S_t$ is mapped onto $Q_{3,3}'$ by the following polynomial automorphism of $\CC^7$
$$
\begin{array}{lll}
z_0^*&=&\displaystyle i\left(z_0-\frac{1}{2}\left(z_1z_6+z_2z_5+z_3z_4\right)-\frac{1}{4}\left(z_4^3+z_5^3+z_6^3+tz_4z_5z_6\right)\right),\\
\vspace{-0.1cm}\\
z_1^*&=&\displaystyle z_1+\frac{3}{2}z_6^2+\frac{t}{2}z_4z_5,\\
\vspace{-0.1cm}\\
z_2^*&=&\displaystyle z_2+\frac{3}{2}z_5^2+\frac{t}{2}z_4z_6,\\
\vspace{-0.1cm}\\
z_3^*&=&\displaystyle z_3+\frac{3}{2}z_4^2+\frac{t}{2}z_5z_6,\\
\vspace{-0.1cm}\\
z_j^*&=&z_j,\quad j=4,5,6.
\end{array}
$$ 

Next, we observe that the base $B_t$ of every hypersurface $S_t$ is affinely homogeneous. Indeed, let us apply to $B_t$ the translation $\tau_q:X\ra X-q$ for $q\in B_t$. Absorbing linear terms into $x_0$, we turn the equation of $\tau_q(B_t)$ into
$$
\begin{array}{l}
x_0={\mathcal Q}_{2}+L_1(x_4,x_5,x_6)x_4+L_2(x_4,x_5,x_6)x_5+L_3(x_4,x_5,x_6)x_6+{\mathcal Q}_3,
\end{array}
$$
where ${\mathcal Q}_j$, $j=2,3$, are the terms of order $j$ in the right-hand side of equation (\ref{hyperst}), and $L_1$, $L_2$, $L_3$ are linear functions. Absorbing $L_1$ into $x_3$, $L_2$ into $x_2$, $L_3$ into $x_1$, we obtain the equation of $B_t$. This proves that $B_t$ is affinely homogeneous. Thus $B_{t_1}$ and $B_{t_2}$ are affinely equivalent if and only if they are linearly equivalent.       

Further, it is straightforward to see that for every $t$ the cubic terms in the equation of $B_t$ are trace-free, where the trace is calculated with respect to the non-degenerate quadratic form ${\mathcal Q}_2$. We now argue as in Section \ref{prooftheorem1} and use Proposition 1 of \cite{EE} to conclude that two hypersurfaces $S_{t_1}$, $S_{t_2}$ are affinely non-equivalent if the corresponding cubics $c_{t_1}$, $c_{t_2}$ (defined in  (\ref{cubics})) are linearly non-equivalent over $\RR$.

To see when two cubics are linearly non-equivalent, we think of them as functions of three complex variables and find sufficient conditions for their zero loci, viewed as curves in $\CC\PP^2$, to be projectively non-equivalent. Namely, for every $t\in\RR$ define
$$
Z_t:=\left\{(w_1:w_2:w_3)\in\CC\PP^2: c_t(w_1,w_2,w_3)=0\right\}.
$$
If $t\ne -3$ the set $Z_t$ is a non-singular elliptic curve, whereas for $t=-3$ it has singularities at the points $(1:q:q^2)$, with $q^3=1$. We only consider small values of $t$ and show that for $t$ lying in some interval around 0 the curves $Z_t$ are pairwise projectively non-equivalent, which implies that for such $t$ the cubics $c_t$ are pairwise linearly non-equivalent, as required.

The projective equivalence class of an elliptic curve (which coincides with its biholomorphic equivalence class) is completely determined by the value of the $j$-invariant for the curve (see e.g. \cite{Ta} or \cite{K}, pp. 56--67). Saito's calculation for simple elliptic singularities of type $\tilde E_6$ in \cite{S2} gives a formula for the value of the $j$-invariant of the curve $Z_t$ (see also a correction made in \cite{CSY}, \cite{E}). Furthermore, it is explained in \cite{E} (see also \cite{CSY}) how one can recover the value of the $j$-invariant directly from the corresponding moduli algebra ${\mathcal A}$. For the completeness of our exposition we compute the value of the $j$-invariant of $Z_t$ below.

In order to apply the well-known formulas for computing the value of the $j$-invariant, we first reduce the equation $Z_t$ to the Weierstrass form (see e.g. formula (1) of \cite{Ta}). [Note that the argument given in \cite{E} does not require such a reduction.]  We perform the following transformation
$$
\begin{array}{lll}
w_1^*&=&\hspace{0.35cm}Cw_3,\\
\vspace{-0.3cm}\\
w_2^*&=&-w_2\displaystyle+\frac{t}{3}w_3,\\
\vspace{-0.3cm}\\
w_3^*&=&\hspace{0.35cm}w_1+w_2-\displaystyle\frac{t}{3}w_3,
\end{array}
$$
where
$$
\displaystyle C:=-\sqrt[3]{\frac{1}{3}\left(\frac{t^3}{27}+1\right)}.
$$
This transformation takes $Z_t$ into the elliptic curve
$$
\displaystyle w_2^2w_3-\frac{t}{3C}w_1w_2w_3+w_2w_3^2=w_1^3-\frac{t^2}{9C^2}w_1^2w_3-\frac{1}{3}w_3^3.
$$
One can now compute the value of the $j$-invariant for this curve as shown in \cite{Ta} (see p. 180). In the notation of \cite{Ta} we have
$$
\displaystyle a_1=-\frac{t}{3C},\,a_2=-\frac{t^2}{9C^2},\,a_3=1,\,a_4=0,\,a_6=-\frac{1}{3},
$$
hence 
$$
\displaystyle c_4=\frac{t}{9C^4}(t^3+72C^3),\,\Delta=\frac{1}{C^3}.
$$
Therefore
$$ 
\displaystyle j=-t^3\frac{(t^3-216)^3}{(t^3+27)^3}.
$$

Setting $s:=t^3$ we obtain
$$
j=\displaystyle\Phi(s):=-s\frac{(s-216)^3}{(s+27)^3}.
$$
The function $\Phi$ is strictly increasing for small $|s|$ since $\Phi'(0)>0$. This shows that the values of the $j$-invariant for the elliptic curves $Z_t$ are pairwise distinct for $t$ lying in a sufficiently small interval around 0.

The proof of the theorem is complete.\qed

{\obeylines
Department of Mathematics
The Australian National University
Canberra, ACT 0200
AUSTRALIA
E-mail: alexander.isaev@maths.anu.edu.au
}

\end{document}